\documentclass{iopjournal}

\usepackage{amsmath,amssymb,amsfonts}
\usepackage{algorithm,setspace}
\usepackage{algpseudocode}
\usepackage{graphicx}

\DeclareMathOperator{\pos}{pos}
\DeclareMathOperator{\prox}{prox}
\DeclareMathOperator*{\argmin}{arg\,min}

\begin{document}

\articletype{Paper} %	 e.g. Paper, Letter, Topical Review...

\title{Efficient primal-dual algorithm for imaging applications with matrix stacking, applied to DBT image reconstruction}

\author{Emil Y. Sidky$^1$$^*$\orcid{0000-0002-6951-2456}, John Paul Phillips$^1$\orcid{0000-0002-3437-425X},
Zheng Zhang$^1$\orcid{0000-0002-6975-6274}, Dan Xia$^1$\orcid{0000-0001-8567-1038},
Ingrid S. Reiser$^1$\orcid{0000-0002-2047-2190}, and Xiaochuan Pan$^1$$^*$\orcid{0000-0002-3074-9771}}

\affil{$^1$Department of Radiology, The University of Chicago, Chicago, Illinois, USA}

\affil{$^*$Author to whom any correspondence should be addressed.}

\email{sidky@uchicago.edu, xpan@uchicago.edu}


\begin{abstract}
The primal-dual hybrid gradient (PDHG) algorithm for solving convex optimization problems that arise in
tomographic imaging is revisited. In particular, simplification of the selection of step-size parameters
is developed for optimization problems with multiple terms, each containing a linear transform subject to splitting.
This simplification maintains algorithm efficiency while avoiding massive grid searches for the optimal
parameter settings. The PDHG framework is demonstrated on an image reconstruction problem for wide-angle
digital breast tomosythesis (DBT); use of the proposed optimization problem is enabled by the framework and it
is demonstrated to have some advantage in quantitative accuracy of the reconstructed volume and
in improving DBT depth resolution.
\end{abstract}

\section{Introduction}
\label{sec:intro}

The PDHG algorithm of Chambolle and Pock \cite{chambolle2011first} has proven to be quite useful for solving non-smooth convex
optimization problems of interest for imaging applications.  In particular, this algorithm has been effective at solving
the large-scale optimization problems that arise in tomographic image reconstruction \cite{sidky2012convex}. 
Originally, the PDHG algorithm was used to address two-term non-smooth problems where one term enforced data-fidelity and
the other term  either penalized or constrained a the image total variation (TV), with the second term being convex but
non-smooth.
More recently, there are
examples of non-smooth optimization applied to tomographic image reconstruction that motivate developing
a system to apply PDHG to convex optimization with multiple terms each involving a different linear transform.
For example, in list-mode time-of-flight (TOF) positron emission tomography (PET), one optimization problem of interest involves
three linear transforms:
a gradient matrix in a TV regularization term plus two
projection matrices, one for activity projection along a line-of-response and the other for the weighted projection for
evaluating a TOF bin estimate \cite{zhang2018optimization}. For digital breast tomosynthesis (DBT) image reconstuction
with directional TV regularization, one optimization of interest involves five linear transforms:
three differential matrices for each of the three spatial directions, an identity matrix for the dualized l1-norm penalty,
and the DBT projection matrix \cite{sidky2025accurate}. With increasing numbers of linear transforms, optimizing
step-size parameters for PDHG algorithm efficiency becomes more and more difficult. This work addresses this issue
directly and it proposes an alternate step-size parameterization that simplifies the parameter value selection.

The paper is organized as follows. Sec.~\ref{sec:PDHG} develops the modified step-size parameterization for the general
PDHG algorithm and suggests a particular form for X-ray based tomographic image reconstruction.
Sec.~\ref{sec:2DDBT} suggests a new form of sparsity-regularization using directional TV (DTV) penalties
for DBT image reconstruction in a 2D setting, and the proposed form of PDHG is used to
solve the corresponding optimization problem. The effect of the proposed step-size
parameters is demonstrated with this system. Finally, in Sec.~\ref{sec:3DDBT} 3D DBT image reconstruction is
addressed for a realistic simulation of wide-angle DBT.

\section{PDHG algorithm for multiple-term non-smooth convex optimization}
\label{sec:PDHG}

\begin{algorithm}
\hrulefill
\begin{algorithmic}[1]
\State $x^{(k+1)} \leftarrow \prox[\tau G] \left( x^{(k)} - \tau K^\top \lambda^{(k)} \right) $
\State $\tilde{x} = 2 x^{(k+1)} - x^{(k)} $
\State $\lambda^{(k+1)} \leftarrow \prox[\sigma F^\star] \left( \lambda^{(k)} + \sigma K \tilde{x} \right) $
\end{algorithmic}
\hrulefill
\caption{PDHG algorithm update steps. Index $k$ indicates the iteration number; $x$ and $y$ are initialized to zero
$x^{(0)} =0$ and $y^{(0)} = 0$; and the parameters $\tau$ and $\sigma$ are the primal and dual step-sizes, respectively.
}
\label{alg:CP1}
\end{algorithm}

The PDHG algorithm solves optimization problems of the following
general form
\begin{equation}
\label{cpopt1}
\min_x \left\{ F(Kx) +G(x) \right\},
\end{equation}
where $F$ and $G$ convex functions, possibly non-smooth; $K$ is a linear transform; and $x$ is the unknown vector, which
for imaging applications, represents the array of pixel/voxel gray values.
The algorithm for solving Eq.~(\ref{cpopt1}) is given in Algorithm~\ref{alg:CP1}.
To derive specific instances of this algorithm, the convex conjugate of $F$ is needed
\begin{equation*}
F^\star (\lambda) = \max_{y} \left\{ \lambda^\top y - F(y) \right\},
\end{equation*}
and the proximal optimizations need to be derived
\begin{align*}
\prox[\tau G](x) & = \argmin_{x^\prime} \left\{ \tau G(x^\prime) + \tfrac{1}{2} \| x^\prime - x \|^2 \right\}, \\
\prox[\sigma F^\star](\lambda) & = \argmin_{\lambda^\prime} \left\{ \sigma F^\star(\lambda^\prime) +
\tfrac{1}{2} \| \lambda^\prime - \lambda \|^2 \right\}.
\end{align*}
The PDHG algorithm is proven to converge as long as the product of the primal and dual step-sizes obeys the condition
$ \sigma \tau < 1/ \|K\|^2_2$, where $\|K\|_2$ is the largest singular value of the matrix $K$, which can be computed
by the power method.
In practice, setting the step-sizes so that this condition is extended to an equality
\begin{equation*}
\sigma \tau = \frac{1}{\|K\|^2_2}
\end{equation*}
leads to convergent PDHG iteration. With this condition, the step-sizes are set to
\begin{equation*}
\sigma = \frac{\beta}{\|K\|_2}, \text{ and } \tau = \frac{1}{\beta \|K\|_2}
\end{equation*}
where $\beta$ is a parameter that needs to be tuned to optimize efficiency.

In imaging, the PDHG algorithm was motivated by the desire to investigate sparsity-regularization with the total variation (TV)
norm. For CT image reconstruction, a common optimization problem studied was TV penalized least-squares with a non-negativity
constraint
\begin{equation}
\label{tvpen}
\min_f \left\{ \tfrac{1}{2} \|Xf - g\|^2_2 + \alpha \| \nabla f \|_1 \text{ subject to } f \ge 0 \right\},
\end{equation}
where $f$ is the image vector; $g$ is sinogram vector; $X$ is a matrix encoding X-ray projection; $\nabla$ is a matrix
encoding the discrete spatial gradient. In this optimization problem, there are two linear transforms
$X$ and $\nabla$. Accordingly, it
is not immediately obvious how to fit this optimization in the form of Eq.~(\ref{cpopt1}).

\subsection{The matrix stacking strategy}
\label{sec:stacking}

One strategy is to stack the matrices and make the following associations
\begin{align*}
x &= f, \text{  }
y =
\left[
\begin{array}{c}
y_\text{proj} \\
y_\text{grad}
\end{array} \right], \\
F(y) &= \frac{1}{2} \|y_\text{proj} - g\|^2 + \frac{\alpha}{\nu} \| y_\text{grad}\|_1, \text{  }
G(x) = \delta(x \,|\, x \ge 0),\\
K &= \left[
\begin{array}{c}
X \\
\nu \nabla
\end{array} \right],
\end{align*}
where the row vectors of $K$ are image vectors and the column vectors are sinograms concatenated with image gradients;
$F$ is composed of two separate terms; $G$ encodes the non-negativity constraint with an indicator function; and
the matrix $K$ is the concatenation of $X$ and $\nabla$. Note that there is another parameter introduced in implementing
this concatenation; the parameter $\nu$ multiplies $\nabla$ in $K$ and it is divided out in $F$. Thus $\nu$ does not
alter the optimization problem, but it will impact algorithm convergence efficiency. As a result there are now two
parameters $\beta$ and $\nu$ the need to be tuned to optimize algorithm convergence. Note that the role of the penalty parameter
$\alpha$ is different from $\beta$ and $\nu$; changing $\alpha$ impacts both the solution of Eq.~(\ref{tvpen})
and algorithm convergence. To thoroughly optimize algorithm efficiency, a 2D parameter search should be performed
over $\beta$ and $\nu$.

For the general case of $N$ terms with different linear transforms, the optimization of interest becomes
\begin{equation}
\label{cpopt2}
\min_x \left\{ \sum_{i=1}^N F_i(K_i x) +G(x) \right\}.
\end{equation}
This optimization can be put in the form of Eq.~(\ref{cpopt1}) with the following associations
\begin{equation*}
y =
\left[
\begin{array}{c}
y_1 \\
y_2 \\
\vdots \\
y_N
\end{array} \right], \; \; \;
\hat{F}(y)=  \sum_{i=1}^N \hat{F}_i(y_i) =  \sum_{i=1}^N F_i(y_i/\nu_i), \; \; \;
K = \left[
\begin{array}{c}
\nu_1 K_1/ \|K_1\|_2 \\
\nu_2 K_2/ \|K_2\|_2 \\
\vdots \\
\nu_N K_N/ \|K_N\|_2
\end{array} \right].
\end{equation*}
For $N$ terms, the number of algorithm parameters is also $N$. As written, there are $N+1$ because
there is the step-size ratio parameter $\beta$ in addition to the $\nu_i$ parameters, but
the parameterization is redundant because the overall scale of the $\nu_i$'s can be absorbed in $\beta$.

It is clear that parameter tuning becomes more laborious as the number of terms increase. One strategy
that has been effective with TV sparsity-regularization is to balance the matrix norms; i.e. set
$\nu = \|X\|_2/\|\nabla\|_2$. This strategy can be extended to the $N$-transform case by
setting all $\nu_i$ to 1,
leaving only the step-size ratio $\beta$ to be tuned.
Though convenient, there is no physical motivation for this and this strategy may be far from optimal in terms of maximizing PDHG efficiency.
Thus investigating other ways to manage multiple linear transforms with PDHG may be beneficial.

\subsection{The generalized proximal point algorithm and step matrices}

\begin{algorithm}
\hrulefill
\begin{algorithmic}[1]
\State $x^{(k+1)} \leftarrow (I + T \partial G)^{-1} \left( x^{(k)} - T K^\top \lambda^{(k)} \right) $
\State $\tilde{x} = 2 x^{(k+1)} - x^{(k)} $
\State $\lambda^{(k+1)} \leftarrow (I + \Sigma \partial F^\star)^{-1} \left( \lambda^{(k)} + \Sigma K \tilde{x} \right) $
\end{algorithmic}
\hrulefill
\caption{Same PDHG update steps as Algorithm~\ref{alg:CP1} except that the scalar step-size
parameters $\sigma$ and $\tau$ are extended to step matrices $\Sigma$ and $T$. The matrix $I$ is the
identity, and $\partial G$ is the subgradient of $G$. The operator $(I + T \partial G)^{-1}$ is
referred to as the resolvent, which is closely linked to the $\prox$ operator.
}
\label{alg:CP2}
\end{algorithm}

An important development for the PDHG algorithm was developed in He and Yuan \cite{he2012convergence},
where it was observed that the PDHG iteration can be put in the form of a generalized proximal point algorithm
allowing the generalization of the step-size parameters $\sigma$ and $\tau$ to symmetric positive definite
step matrices $\Sigma$ and $T$, yielding the PDHG iteration shown in Algorithm~\ref{alg:CP2}.
The resolvents in lines 1 and 3 of this algorithm arise from generalizing the $\prox$ optimization.
With the step matrix $T$, for example, this optimization becomes
\begin{equation*}
\argmin_{x^\prime} \left\{ G(x^\prime) + \tfrac{1}{2} (x^\prime - x)^\top T^{-1}(x^\prime - x)  \right\}.
\end{equation*}
Solving this optimization is equivalent to applying the resolvent $(I + T \partial G)^{-1}$, when $T$
is positive definite.
% In this work, we will consider $T$ that is positive semi-definite in which
%case the resolvent is still well-defined but the corresponding optimization may not be because $T$
%is not invertible if it has $0$ eigenvalues.
%In addition, this work proposed predictor-corrector extensions to the basic
%PDHG algorithm that can improve algorithm efficiency.
Pock and Chambolle \cite{pock2011diagonal} showed that the iteration in Algoritm~\ref{alg:CP2}
converges if 
\begin{equation}
\label{matrixcond}
\|\Sigma^{1/2} K T^{1/2}\|^2_2 = \|T^{1/2} K^\top \Sigma K T^{1/2}\|_2 < 1,
\end{equation}
and in practice this condition can be extended to an equality. 
In their work, they use the generalization to matrix steps to derive diagonal preconditioners, which
are less computationally demanding than using the power method to arrive at the scalar step-sizes in
the original PDHG formulation. The proposed step matrices do allow additional flexibility
in the matrix stacking strategy that can be exploited for efficiency and intuitiveness of determining
algorithm parameters.

\begin{algorithm}
\hrulefill
\begin{algorithmic}[1]
\State $x^{(k+1)} \leftarrow (I + \tau \partial G)^{-1} \left( x^{(k)} - \tau \left(\sum_i \nu_i \hat{K}_i^\top \lambda_i^{(k)}\right) \right) $
\State $\tilde{x} = 2 x^{(k+1)} - x^{(k)} $
\State $\lambda_i^{(k+1)} \leftarrow (I + \sigma_i \partial \hat{F}_i^\star)^{-1} \left( \lambda_i^{(k)} + \sigma_i \nu_i \hat{K}_i \tilde{x} \right)  \;\;\; \forall i$
\end{algorithmic}
\hrulefill
\caption{PDHG update steps derived from Algorithm~\ref{alg:CP2} using the step-size matrices stated
in Eq.~\ref{blocksteps}. The normalized matrices $\hat{K}_i = K_i/\|K_i\|_2$ are used for conciseness.
}
\label{alg:CP3}
\end{algorithm}

\subsection{The matrix stacking strategy revisited}
\label{sec:newstacking}

With the generalization to step matrices, it becomes possible to adjust step-size parameters for each
of the individual matrices comprising $K$ by writing $T$ as a scalar times
the identity matrix and $\Sigma$ in the following diagonal form
\begin{equation}
\label{blocksteps}
T = \tau I, \; \; \;
\Sigma = 
\left( 
\begin{array}{cccc}
\sigma_1 I & \cdot     & \dots  & \cdot      \\
\cdot      &\sigma_2 I & \dots  & \cdot      \\
\vdots     & \vdots    & \ddots & \vdots     \\
\cdot      & \cdot     & \dots  & \sigma_N I
\end{array}
\right).
\end{equation}
Together with the matrix scaling parameters $\nu_i$ the condition stated in Eq.~(\ref{matrixcond})
becomes
\begin{equation}
\label{stepcondition}
\tau \left\| \sum_i \sigma_i \nu^2_i \hat{K}^\top_i \hat{K}_i \right \|_2 = 1,
\end{equation}
where we use the normalized block matrices $\hat{K}_i = K_i/\|K_i\|_2$.

The generic PDHG optimization for stacked matrices is stated in Eq.~(\ref{cpopt2}), and
the PHDG update steps are written down in Algorithm~\ref{alg:CP3}. Note that the coefficient
of the linear term, $K_i^\top y_i^{(k)}$, in the primal update is $\tau \nu_i$ and
the coefficient of the linear term, $K_i\tilde{x}$, in the dual update is $\sigma_i \nu_i$.
Thus by adjusting $\nu_i$ and $\sigma_i$ appropriately we can control the overall relative step-size
and step-size ratio corresponding to each $K_i$ block. 
We define the step weights, and ratios, as the product, and ratios, of these two coefficients
\begin{align}
w_i = \tau \sigma_i \nu_i^2, \label{wdef} \\
r_i = \sigma_i/\tau \label{rdef}
\end{align}
because the parameters $w_i$ and $r_i$ lend themselves better to heuristic arguments for their
determination. The weights $w_i$ are represented as an overall weight $w$ multiplied by
a relative weight $\hat{w}_i$, where $w$ is determined by Eqs.~(\ref{stepcondition}) and
(\ref{wdef})
\begin{equation}
\label{wmag}
w = \left\| \sum_i \hat{w}_i \hat{K}^\top_i \hat{K}_i \right\|^{-1}_2.
\end{equation}
Eliminating $\sigma_i$ from Eq.~(\ref{wdef}) using Eq.~(\ref{rdef}) yields
\begin{equation}
\label{taunu}
w_i = w \, \hat{w}_i =  r_i \tau^2 \nu_i^2,
\end{equation}
and in this equation $w_i$ and $r_i$ are known, but there are two
unknowns $\tau$ and $\nu_i$. To pin down $\tau$, we note that leaving
all $\nu_i$ free to vary amounts to an over-parameterization because
an overall scale change, $k$, of the $\nu_i$ can be compensated for
by a scale change of $1/\sqrt{k}$ in $\sigma_i$ and $\tau$. Therefore, we
fix the scale of the $\nu_i$ by selecting
\begin{equation}
\label{nu1}
\nu_1 = 1,
\end{equation}
allowing the determination of $\tau$ using Eq.~(\ref{taunu}) for $i=1$
\begin{equation}
\label{taudef}
\tau = \sqrt{ w \, \hat{w}_1 /r_1}.
\end{equation}
In turn, knowledge of $\tau$ yields the $\sigma_i$ parameters
\begin{equation}
\label{sigdef}
\sigma_i = r_i \tau,
\end{equation}
and with Eq.~(\ref{taunu}), it also yields the $\nu_i$ parameters
\begin{equation}
\label{nudef}
\nu_i = (1/\tau)\sqrt{ w \, \hat{w}_i / r_i }.
\end{equation}
Thus all of the step parameters for implementing Algorithm~\ref{alg:CP3} are determined
from the desired relative step weights, $\hat{w}_i$, and ratios, $r_i$.

\begin{algorithm}
\hrulefill
\begin{algorithmic}[1]
\State $\hat{w}_1 = \gamma; \; \; \hat{w}_i = 1 \; \;\; \forall i \ge 2$
\State $r_i = \beta \; \; \; \forall i$
\State Step parameters $w$, $\tau$, $\sigma_i$, and $\nu_i$ are obtained from Eqs.~(\ref{wmag}) through (\ref{nudef}).
\For {$k \leftarrow 1$, $N_\text{iter}$}
   \State $\hat{x} \leftarrow (I + \tau \partial G)^{-1} \left( x^{(k)} - \tau \left(\sum_i \nu_i \hat{K}_i^\top \lambda_i^{(k)}\right) \right) $
   \State $\tilde{x} = 2 \hat{x} - x^{(k)} $
   \State $\hat{\lambda}_i \leftarrow (I + \sigma_i \partial \hat{F}_i^\star)^{-1} \left( \lambda_i^{(k)} + \sigma_i \nu_i \hat{K}_i \tilde{x} \right)  \;\;\; \forall i$
   \State $x^{(k+1)} \leftarrow x^{(k)} + \rho \left( \hat{x} - x^{(k)} \right)$
   \State $\lambda_i^{(k+1)} \leftarrow \lambda_i^{(k)} + \rho \left( \hat{\lambda}_i - \lambda_i^{(k)} \right)  \;\;\; \forall i$
\EndFor
\end{algorithmic}
\hrulefill
\caption{PDHG update steps proposed for use on X-ray tomography imaging models. The last two lines
originate from a predictor-corrector algorithm explained in Ref. \cite{he2012convergence}. The first linear
transform block, $K_1$, is based on the X-ray transform $X$. The algorithm parameters, $\beta \in (0,\infty)$, $\gamma \in [1,\infty)$,
and $\rho \in [1,2)$ need to be tuned and serve as inputs to the pseudocode.
}
\label{alg:CP4}
\end{algorithm}

\subsection{Step parameter heuristic for X-ray tomographic image reconstruction}
\label{sec:steps}

For X-ray tomographic image reconstruction, the dominant term in the imaging
model optimization problem, which guides the solution and limits
convergence, is the data-discrepancy term involving the projection matrix $X$.
Accordingly, the proposal is to weight the step-sizes corresponding to this linear transform block
more than the linear transforms corresponding to the various regularization terms.
Specifically, we assign $K_1 = X$, $\hat{w}_1 = \gamma$, and $\hat{w}_i =1$ for $i \ge 2$.
Next, we assign all step-size ratios corresponding to each of the matrix blocks to the same parameter $r_i = \beta$
for all $i$. With this system, only two step-size parameters, $\gamma$ and $\beta$, need to be searched.
This simplifications in step-size parameterization
enables the use of acceleration techniques, which may introduce additional parameters, for improving convergence rate.
In particular, we consider a predictor-corrector
extension to the PDHG algorithm proposed by He and Yuan \cite{he2012convergence}.

Putting it all together, the PDHG algorithmic template, which we work with here, is stated
in Algorithm~\ref{alg:CP4}. There are three step-size parameters to tune. The range for $\gamma$
is fairly restricted, the useful search interval is expected to be $[1,10]$. The setting $\gamma=1$
weights the X-ray transform-based matrix block, $K_1$, on an equal footing with the other blocks of $K$.
Setting $\gamma > 1$, apportions more weighting to $K_1$, but there is limiting returns on increasing $\gamma$.
Likewise, the useful range for the predictor-corrector parameter is $\rho \in [1,2)$. The step-size ratio 
parameter, $\beta$, on the other hand, is the most difficult parameter to optimize because there does not
appear to be any rhyme or reason on how to select it. Accordingly, $\beta$ should be searched on a logarithmic
grid covering many orders of magnitude.

\section{Limited angular-range image reconstruction with directional TV penalties in 2D}
\label{sec:2DDBT}

\begin{figure}[t]
 \centering
        \includegraphics[width=0.9\textwidth]{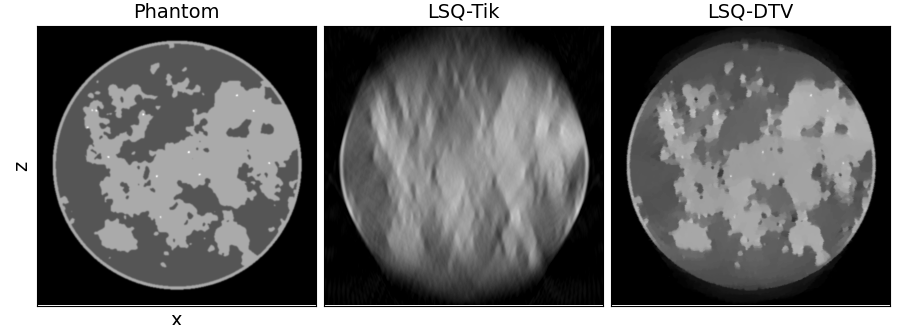}
\caption{(Left) 512x512 test phantom for a 2D breast slice simulation. (Middle) reconstruction of limited-angular
range data for an arc of 50$^\circ$ with Tikhonov regularized least-squares.
(Right) reconstruction with using gradient-sparsity regularization based on directional TV.
The depth direction, labelled $z$, is vertical and the in-plane direction, labelled $x$, is horizontal.
The scanning arc, located above the phantom is centered on the vertical mid-line of the shown images.}
\label{fig:recon2D}
\end{figure}

To demonstrate the utility of the matrix stacking approach proposed in Sec.~\ref{sec:steps}, we apply it to an optimization problem
based on one proposed for tomographic image reconstruction for DBT from Sidky {\it et al.}\cite{sidky2025accurate}.
In that work, data-discrepancy constrained, sparsity-regularization
in the following form was developed
\begin{equation}
\label{DTVopt}
\argmin_f \left\{ \alpha_x \| \partial_x f \|_1 + \alpha_z \| \partial_z f \|_1 + \alpha_1 \|I f \|_1
\text{ such that } \|R[c](Xf-g)\|_2 \le \epsilon \cdot \sqrt{\text{size}(g)} \text{ and } f \ge 0 \right\},
\end{equation}
where the $\ell_1$-norm terms encourage sparsity in the directional gradients of the image and in the image
itself. The $x$- and $z-$directions are the DBT in-plane and depth directions, see Fig.~\ref{fig:recon2D}.
The matrix operators are $\partial_x$, $\partial_z$, $I$, and $X^\prime=R[c]X$; these matrices
represent, respectively, finite differencing of the image in the $x-$ and $z$-directions, the identity matrix,
and the composition of $R[c]$  with the projection matrix, where $R[c]$ encodes filtering with the square
root of the Hanning filter with cutoff parameter $c$. The heuristic for filtering with the square root of the Hanning filter
stems from writing down a gradient descent step for the least-squares data discrepancy term $\tfrac{1}{2}\|R[c](Xf-g)\|^2_2$
\begin{align*}
f^{(k+1)} & = f^{(k)} - \alpha  X^\top \, (R^\top [c] R[c]) \, (Xf^{(k)} - g), \\
f^{(1)} & = \alpha X^\top \, (R^\top [c] R[c]) \, g \;\;\; \text{if} \;\;\; f^{(0)} =0,
\end{align*}
where the first iteration is the FBP image if the the initial image is zero, and the filter is $R^\top [c] R[c]$;
for symmetric filters $R^\top [c]$ = $R[c]$. 
We write the $\ell_1$ penalty term explicitly with the $I$ matrix because we wish to include this term
with the $F(\cdot)$ function instead of the $G(\cdot)$ function of the generic optimization that PDHG solves.
The parameters $\alpha_x$, $\alpha_z$, and $\alpha_1$ are penalty weights, and $\epsilon$ is the data-discrepancy
constraint parameters. The latter parameter is scaled by the square root of the number of elements $\text{size}(g)$
in the data vector so that $\epsilon$ represents a root-mean-square-error measure.

\begin{figure}
 \centering
        \includegraphics[width=0.5\textwidth]{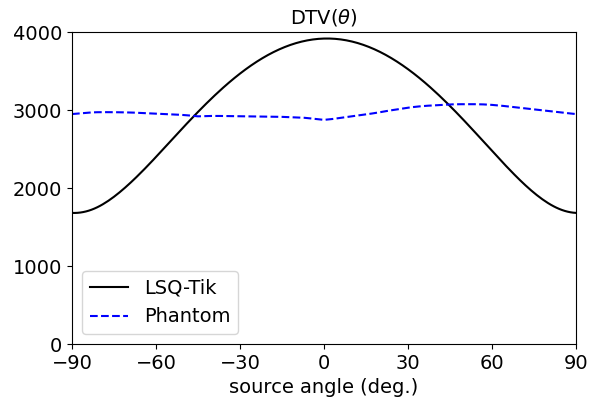}
 \caption{DTV$(\theta)$ of test phantom and LSQ-Tik reconstructed image, shown in Fig.~\ref{fig:recon2D}.}
\label{fig:DTVangle}
\end{figure}

In this work, we modify the directional gradient terms in a way that is more effective than the proposed
optimization in Eq.~\ref{DTVopt} and that allows for intuitive determination of the directional penalty parameters.
In Fig.~\ref{fig:recon2D}, results for a 2D breast DBT simulation are shown where the scanning arc covers a limited
angular range of 50$^\circ$ with 25 projections. The simulated detector is comprised of 1024 detector pixels
arranged on a linear array and the detector rotates opposite the X-ray source as in a CT scan.
The test object is shown in Fig.~\ref{fig:recon2D} and from this discrete image, ideal projection data are generated
with no noise added.  The middle panel of Fig.~\ref{fig:recon2D} shows the image reconstruction result for Tikhonov
regularized least-squares (LSQ-Tik).\cite{rose2017} The Tikhonov regularization parameter is set so that the mean final data-discepancy
is 0.1\%. As expected, there is significant depth-blurring in the LSQ-Tik reconstruction; this image is useful
in that it provides guidance on how to design directional total-variation (DTV) penalties. We define the
directional derivative matrix as
\begin{equation}
\label{dgrad}
\partial_{\theta} = \cos \theta \partial_x + \sin \theta \partial_z,
\end{equation}
where the derivative is designed to be taken in the direction perpendicular to the source-center-of-rotation line and parallel
to the source arc. The angle $\theta$ is referenced to the mid-point of the source arc, hence $\theta = 0$ results in
differentiation along the $x$-axis. The corresponding DTV$(\theta)$ term is simply $\| \partial_\theta f \|_1$.
DTV$(\theta)$ for the LSQ-Tik reconstruction is plotted in Fig.~\ref{fig:DTVangle} and referenced to DTV$(\theta)$
of the phantom. The LSQ-Tik DTV$(\theta)$ is greater than that of the truth in the range $ - 50^\circ \le \theta \le 50^\circ $
and only in this range of $\theta$ will a DTV penalty be directly effective. In our previous work the small penalty
on DTV in the $z$-direction was only helpful because it was used in combination with a DTV penalty in the $x$-direction.
Moreover the ratio of the $x$-DTV and $z$-DTV parameters needed to be tuned.
Based on the observation in Fig.~\ref{fig:DTVangle}, the DTV$(\theta)$ should be selected for $\theta$ in the range where
the LSQ-Tik DTV$(\theta)$ is larger than that of the test phantom. Furthermore, to first approximation the corresponding
penalty parameters can all be set to the same value -- simplifying the penalty parameter search.

The modified constrained, DTV-minimization formulation
involves eliminating the $z$-DTV term and replacing it with two additional DTV$(\theta)$
terms corresponding to the endpoints of the scanning arc; i.e. DTV at $\theta = -25^\circ, \, 25^\circ$.
The full imaging model is specified by the following optimization
\begin{multline}
\label{DTVoptNew}
\argmin_f \left\{ \alpha_x \| \partial_x f \|_1 + \alpha_a \| \partial_a f \|_1  + \alpha_b \| \partial_b f \|_1 + \alpha_1 \|I f \|_1 \right. \\
\left.
\text{ such that } \|R[c](Xf-g)\|_2 \le \epsilon \cdot \sqrt{\text{size}(g)} \text{ and } f \ge 0 \right\},
\end{multline}
where we use the labels $a$ and $b$ to refer to the directions $\theta = -25^\circ$ and $25^\circ$, respectively;
and we constrain the sum of the penalty parameters to equal one, because multiplying all $\alpha$ by a constant
does not alter the solution
\begin{equation*}
\alpha_x + \alpha_a + \alpha_b + \alpha_1 =1.
\end{equation*}
This optimization problem highlights the need for a matrix stacking strategy because there are a total of five
matrices in the various terms of Eq.~(\ref{DTVoptNew}).

\begin{algorithm}
\hrulefill
\begin{algorithmic}[1]
\State $\hat{f} \leftarrow  f^{(k)} -
\tau \left( (X^\prime)^\top \lambda_s^{(k)}/L_s + \nu_x (\partial_x)^\top \lambda_x^{(k)}/L_x
 + \nu_a (\partial_a)^\top \lambda_a^{(k)}/L_a + \nu_a (\partial_b)^\top \lambda_b^{(k)}/L_b + \nu_1 \lambda_1^{(k)} \right) $
\vspace{2mm}
\State $\hat{f} \leftarrow \text{pos}(\hat{f})$
\vspace{2mm}
\State $\tilde{f} = 2 \hat{f} - f^{(k)} $
\vspace{2mm}
\State $\hat{\lambda}_s^{(k+1)} \leftarrow \text{max}
\left( \left\|\lambda_s^{(k)} + \sigma_s ( X^\prime \tilde{f}/L_s - g^\prime) \right\|_2 - \sigma_s \epsilon^\prime, 0 \right)
\cdot \text{norm}\left(\lambda_s^{(k)} + \sigma_s (\nu_s X^\prime \tilde{f}/L_s - g^\prime)\right) $% \cdot
%\left\|\lambda_s^{(k)} + \sigma_s ( X^\prime \tilde{f} - g^\prime)\right\|^{-1}_2 $
\vspace{2mm}
\State $\hat{\lambda}_x \leftarrow \alpha^\prime_x \left(\lambda_x^{(k)} + \nu_x \sigma_x \partial_x \tilde{f}/L_x \right)/
\text{max} \left( \alpha^\prime_x, \left|  \lambda_x^{(k)} + \nu_x \sigma_x \partial_x \tilde{f}/L_x \right| \right)$
\vspace{2mm}
\State $\hat{\lambda}_a \leftarrow \alpha^\prime_a \left(\lambda_a^{(k)} + \nu_a \sigma_a \partial_a \tilde{f}/L_a \right)/
\text{max} \left( \alpha^\prime_a, \left|  \lambda_a^{(k)} + \nu_a \sigma_a \partial_a \tilde{f}/L_a \right| \right)$
\vspace{2mm}
\State $\hat{\lambda}_b \leftarrow \alpha^\prime_b \left(\lambda_b^{(k)} + \nu_b \sigma_b \partial_b \tilde{f}/L_b \right)/
\text{max} \left( \alpha^\prime_b, \left|  \lambda_b^{(k)} + \nu_b \sigma_b \partial_b \tilde{f}/L_b \right| \right)$
\vspace{2mm}
\State $\hat{\lambda}_1 \leftarrow \alpha^\prime_1 \left(\lambda_1^{(k)} + \nu_1 \sigma_1 \tilde{f} \right)/
\text{max} \left( \alpha^\prime_1, \left|  \lambda_1^{(k)} + \nu_1 \sigma_1  \tilde{f} \right| \right)$
\vspace{2mm}
\State $f^{(k+1)} \leftarrow f^{(k)} + \rho \left( \hat{f} - f^{(k)} \right)$
\vspace{2mm}
\State $\lambda_s^{(k+1)} \leftarrow \lambda_s^{(k)} + \rho \left( \hat{\lambda}_s - \lambda_s^{(k)} \right)$
\vspace{2mm}
\State $\lambda_x^{(k+1)} \leftarrow \lambda_x^{(k)} + \rho \left( \hat{\lambda}_x - \lambda_x^{(k)} \right)$
\vspace{2mm}
\State $\lambda_a^{(k+1)} \leftarrow \lambda_a^{(k)} + \rho \left( \hat{\lambda}_a - \lambda_a^{(k)} \right)$
\vspace{2mm}
\State $\lambda_b^{(k+1)} \leftarrow \lambda_b^{(k)} + \rho \left( \hat{\lambda}_b - \lambda_b^{(k)} \right)$
\vspace{2mm}
\State $\lambda_1^{(k+1)} \leftarrow \lambda_1^{(k)} + \rho \left( \hat{\lambda}_1 - \lambda_1^{(k)} \right)$
\end{algorithmic}
\hrulefill
\caption{PDHG algorithm update steps for solving Eq.~(\ref{DTVoptNew}).
The step-size parameters $\beta$, $\gamma$, and $\rho$ are tuned.
The $\sigma$ and $\nu$ parameters and $\tau$ are determined as described in the caption of Algorithm~\ref{alg:CP4}.
The stated update steps in this algorithm correspond to lines 5-9 inside the for-loop of Algorithm~\ref{alg:CP4}.
The function $\pos(\cdot)$ thresholds negative values of its argument to zero.
We use $\text{norm}(\cdot)$ to indicate a function that normalizes the argument vector to unit length.
In lines 3-5 the function $\text{max}(\cdot)$ operates component-wise on the second argument 
with the scalar in the first argument.
}
\label{alg:DTV}
\end{algorithm}

\subsection{Matrix stacking approach to constrained, DTV-minimization}

The PDHG algorithm for solving
Eq.~(\ref{DTVoptNew}) is derived. The first step toward this derivation is making the associations necessary
for realizing a PDHG instance
\begin{align}
x &= f, \text{  }
y =
\left[
\begin{array}{c}
y_s \\
y_{x}\\
y_{a}\\
y_{b}\\
y_{1}\\
\end{array} \right], \; \; L_s = \|X^\prime\|_2=\|R[c] X\|_2, \; L_x = \|\partial_x\|_2, \; L_a = \|\partial_a\|_2,
\; L_b = \|\partial_b\|_2 \notag \\
F(y) &=  \delta (y_s \, |\, \|y_s - g^\prime\|_2 \le \epsilon^\prime) + \alpha_x^\prime \| y_{x}\|_1 +
\alpha_a^\prime \| y_{a}\|_1 +\alpha_b^\prime \| y_{b}\|_1 +  \alpha_1^\prime \| y_{1}\|_1, \notag \\
& \text{where }\\
g^\prime = & R[c] g/L_s, \;\; \epsilon^\prime = \epsilon \cdot \sqrt{\text{size}(g)}/L_s, \; \; 
\alpha_x^\prime = \alpha_x L_x/\nu_x, \;\; \alpha_a^\prime = \alpha_a L_a/\nu_a, \;\; \alpha_b^\prime = \alpha_b L_b/\nu_b,
\;\; \alpha_1^\prime = \alpha_1/\nu_1, \notag \\
G(f) &= \delta(f \,|\, f \ge 0), \; \; \;
K = \left[
\begin{array}{c}
X^\prime/L_s \\
\nu_x \partial_x /L_x\\
\nu_a \partial_a /L_a\\
\nu_b \partial_b /L_b\\
\nu_1 I
\end{array} \right], \label{msk}
\end{align}
where no $\nu$ factor multiplies $X^\prime$ by the convention adopted in Eq.~(\ref{nu1}).
Computing the convex conjugate for each term in $F(\cdot)$ yields
\begin{align*}
F^*_s(\lambda_s) &= \epsilon^\prime\|\lambda_s\|_2 + \lambda_s^\top g^\prime, \\
F^*_x(\lambda_x) &= \delta ( \lambda_x \, | \, \|\lambda_x\|_\infty \le \alpha_x^\prime), \\
F^*_a(\lambda_a) &= \delta ( \lambda_a \, | \, \|\lambda_a\|_\infty \le \alpha_a^\prime), \\
F^*_b(\lambda_b) &= \delta ( \lambda_b \, | \, \|\lambda_b\|_\infty \le \alpha_b^\prime), \\
F^*_1(\lambda_1) &= \delta ( \lambda_1 \, | \, \|\lambda_1\|_\infty \le \alpha_1^\prime),
\end{align*}
where
\begin{equation*}
\|v\|_\infty = \text{max} |v|;
\end{equation*}
i.e. the $\ell_\infty$-norm is the magnitude of the largest component of $|v|$.
Because the $F^*$ and $G$ functions are simple, their correspond prox optimizations,
needed to instantiate Algorithm~\ref{alg:CP4}, can be analytically
derived, yielding the PDHG update steps in Algorithm~\ref{alg:DTV}.

\begin{figure}
 \centering
        \includegraphics[width=0.45\textwidth]{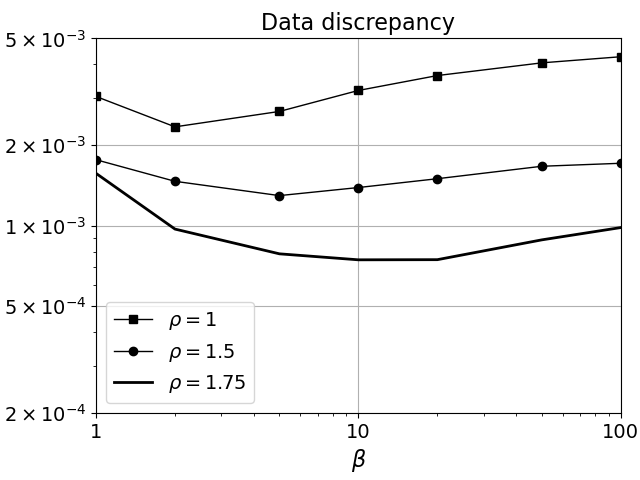}
        \includegraphics[width=0.45\textwidth]{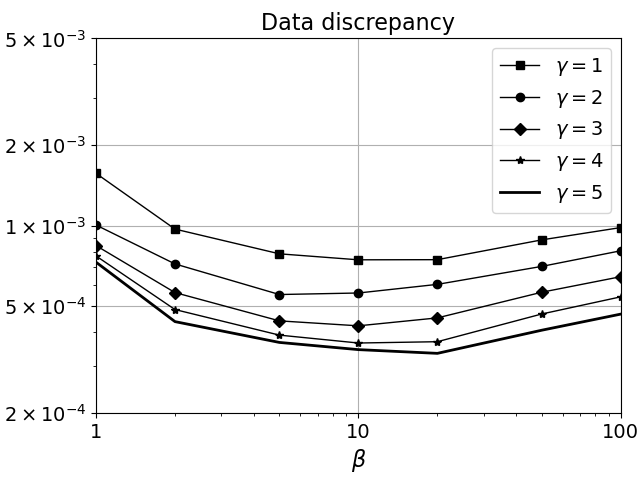}
 \caption{Data residual after 1,000 iterations of Algorithm~\ref{alg:DTV}}
\label{fig:rho}
\end{figure}

\subsection{Step-size parameter tuning for the 2D breast DBT simulation}

For the purpose of demonstrating the new step-size parameterization, we return to the simulation
depicted in Fig.~\ref{fig:recon2D}. Generating ideal, noiseless simulation data, we perform image
reconstruction with LSQ-DTV setting the data error constraint parameter to a small value $\epsilon= 10^{-6}$.
Executing 1,000 iterations, we measure convergence by how small is the resulting data discrepancy.
For all cases, the $\alpha$ penalty parameters are set to equal values, i.e.
\begin{equation*}
\alpha_x = \alpha_a = \alpha_b = \alpha_1 =1/4.
\end{equation*}
The resulting convergence plots as a function of step-size parameters $\beta$, $\rho$, and $\gamma$ are
shown in Fig.~\ref{fig:rho}. The results of the left panel of Fig.~\ref{fig:rho} show the effectiveness
of the He-Yuan predictor-corrector extension to PDHG, and the results of the right panel shows additional
improvement in convergence gained by tuning $\gamma$. The corresponding reconstructed image is shown in the
right panel of Fig.~\ref{fig:recon2D}, where it is seen that sparsity regularization helps to recover
the test phantom.

\section{Application of directional TV penalties to wide-angle DBT image reconstruction}
\label{sec:3DDBT}

Having provided the details for implementing PDHG on a 2D DBT simulation, we are in a position to address
3D DBT for a realistic wide-angle configuration. The imaging model considered here is a simple modification
to Eq.~(\ref{DTVoptNew})
\begin{multline}
\label{DTVopt3D}
\argmin_f \left\{ \alpha_x \| \partial_x f \|_1 + \alpha_y \| \partial_y f \|_1  +
\alpha_a \| \partial_a f \|_1  + \alpha_b \| \partial_b f \|_1 + \alpha_1 \|I f \|_1 \right. \\
\left.
\text{ such that } \|R[c](XG[d]f-g)\|_2 \le \epsilon \cdot \sqrt{\text{size}(g)} \text{ and } f \ge 0 \right\},
\end{multline}
where a DTV$_y$ term has been included that corresponds to the in-plane $y$-direction, which
is perpendicular to the source trajectory; and gaussian blurring operation has been included in data-error
constraint term, which involves blurring the underlying image $f$ with a 3D gaussian of width $d=(d_x,d_y,d_z)$.
Modeling the image as a gaussian-blurred pixel array is particularly useful for low-resolution image
reconstruction to avoid a jagged image representation. The details for deriving the corresponding
PDHG algorithm are not repeated here because the derivation and PDHG instantiation are a straightforward
extension of the presentation in Sec.~\ref{sec:2DDBT}.

\begin{figure}
 \centering
        \includegraphics[width=0.45\textwidth]{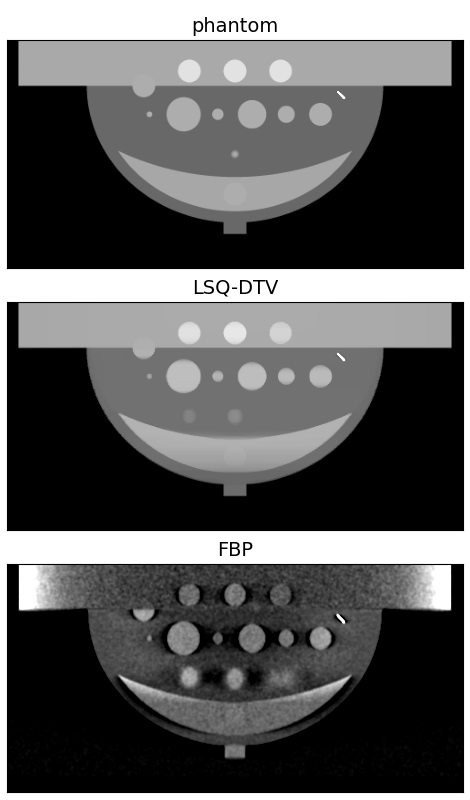}
        \includegraphics[width=0.45\textwidth]{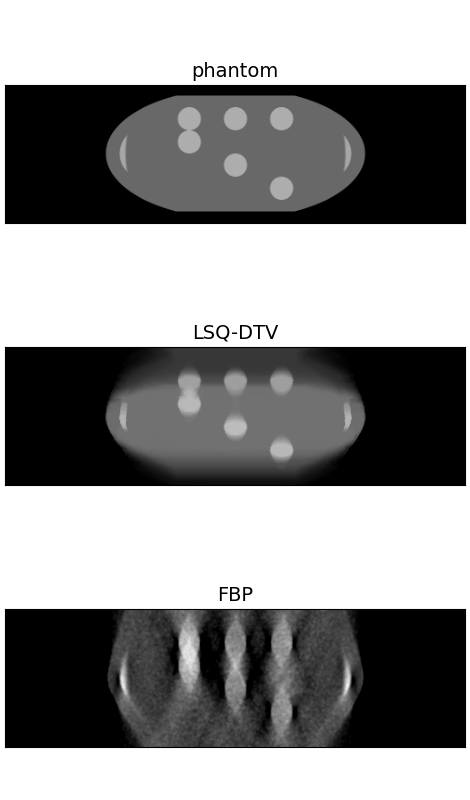}
 \caption{Reconstructed images of the geometric-shape breast phantom: in-plane slices are shown on
 the left and depth-plane slices with $z$-resolution test signals are shown on the right.
 The top images are the phantom slices; the middle row are the LSQ-DTV results; and the bottom images
 show the corresponding FBP reconstruction. The gray scale window for all images is [0,0.6] cm$^{-1}$,
 but the FBP images are scaled so that the phantom objects are clearly visible in this gray scale window.}
\label{fig:phgeom}
\end{figure}

\subsection{Toward low-resolution quantitative DBT image reconstruction}
\label{sec:lowres}

For the first study, we demonstrate the updated PDHG algorithm for LSQ-DTV on a realistic wide-angle
DBT setup with a breast phantom composed of geometric shapes \cite{reiserDBTphantom2}, seen in the top row of Fig.~\ref{fig:phgeom}.
Because the phantom is composed of geometric shapes with analytically computable projections, it is
useful for investigating discretization error as well as providing a test for contrast and depth-resolution.
As discussed in our previous work \cite{sidky2025accurate}, quantitative image reconstruction in DBT may be
possible for simple distributions using voxel sizes much larger than the detector pixel size.
Accordingly, we simulate a set-up
where the detector pixel size is 100 microns squared and the voxel size is 400 microns squared.
With such large voxels, discretization error can be correspondingly large. To control this error, the
gaussian blur kernel in Eq.~(\ref{DTVopt3D}) is assigned a width $d=(400,400,400)$ microns; additionally,
the projections in the data set $g$ are blurred with a gaussian kernel of width $(500,500)$ microns.
The projection blurring is essential for filtering out data components that cannot be represented by the large voxel
image representation. The width of the data blurring kernel can be determined by noting the accessible
data-error parameter, $\epsilon$, values. With the stated set-up, $\epsilon$ can be selected below $0.001$.

The remainder of the simulation parameters are as follows. As quantitative DBT relies on as large of a scanning arc
as possible, we model the existing wide-angle DBT system using 25 projections covering a 50$^\circ$. The detector
array, modeled as 3200$\times$1000 pixels of width 100 microns, is large enough to avoid projection truncation.
For including noise, a Poisson distribution is modeled corresponding to a fluence of 50,000 photons incident on
each detector pixel. Mean sinogram values are computed by analytic line-integration
of the geometric-shapes breast phantom. The image array is composed of 500$\times$250$\times$150 cubic voxels 400 microns.

The parameters of the optimization, except for $\epsilon$, are searched on a coarse grid minimizing the image RMSE.
The resulting $\alpha$ values are chosen as
\begin{equation*}
\alpha_y = \alpha_a = \alpha_b = \alpha_1 =1/9; \; \; \alpha_x = 5/9,
\end{equation*}
respecting the constraint that they sum to 1.
The data-error parameter $\epsilon$ has a strong impact on the image reconstruction results and must be tuned
on a fine grid. For the present results it is selected to be $\epsilon = 0.0011$.
The algorithm parameters are selected as $\beta = 50$, $\gamma = 5$, and $\rho = 1.75$. The parameters $\gamma$ and
$\rho$ are chose to be the same as for the 2D simulation in Sec.~\ref{sec:2DDBT} and optimal $\beta$ is searched based
the number of iterations needed to reach the data-error constraint value.

\begin{figure}
 \centering
        \includegraphics[width=0.5\textwidth]{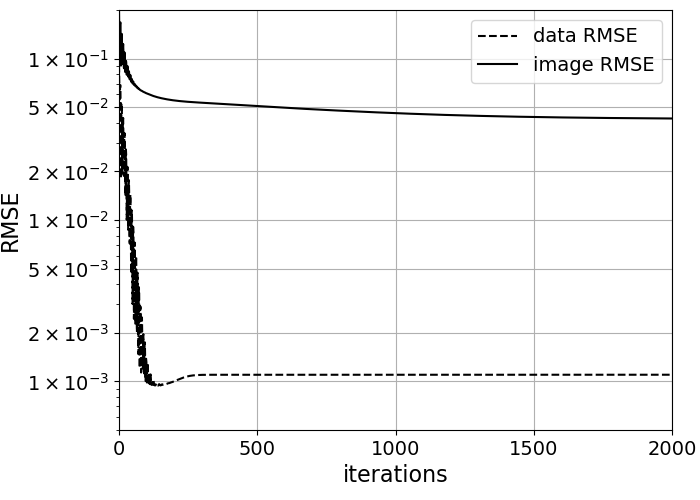}
 \caption{Filtered data and image RMSE corresponding to the LSQ-DTV results shown in Fig.~\ref{fig:phgeom}}
\label{fig:phgeomrmse}
\end{figure}

In-plane and depth slices for the resulting reconstructed volume is shown in Fig.~\ref{fig:phgeom}, where the phantom
and LSQ-DTV result are shown with a consistent gray scale window so that a quantitatve comparison can be made.
Filtered back-projection (FBP) images are shown for reference; the in-plane FBP slice image shows significant contamination
from signals above and below the viewing plane and depth-plane shows the difficulty in resolving the spherical signals
in the depth direction $z$. The FBP images are not regularized so that the noise level can be appreciated. The
LSQ-DTV result shows decent quantitative agreement on the selected in-plane slice, although there are visible differences
from the breast phantom image. The depth-plane slice comparison reveals more obvious discrepancies for the LSQ-DTV
result particularly in the breast background near the top and bottom of the breast. Nevertheless, LSQ-DTV is able
to resolve the circular signals in the $z$-direction with the FBP result indicating the difficulty of this imaging
task.

\begin{figure}
 \centering
        \includegraphics[width=0.5\textwidth]{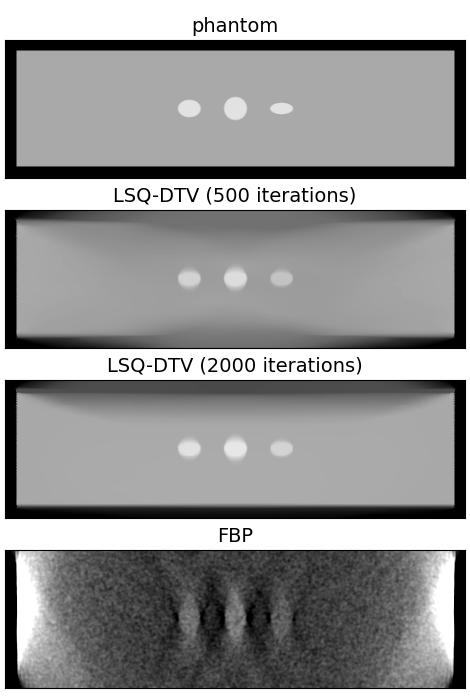}
 \caption{Depth-plane slice images through the three signals in the muscle wall section of the breast phantom.
 The gray scale is [0, 0.6] cm$^{-1}$ for all images and the FBP result is scaled so that the structures are visible
 in this window.}
\label{fig:phgeommw}
\end{figure}

The obtained results with the modified optimization and improved PDHG implementation
shows some improvement compared to a similar study carried out in Sidky {\it et al.} \cite{sidky2025accurate}, although
the previous study did not include noise.
The more relevant result for the present purpose relates to convergence. In Fig.~\ref{fig:phgeomrmse}, the RMSEs
of the filtered data-discrepancy and the image are plotted over 2000 iterations of LSQ-DTV.
The data RMSE constraint value of $\epsilon=0.0011$ is reached at about 200 iterations, but the image RMSE decreases
slowly over the shown 2000 iterations. Also, this gradual improvement in the image RMSE corresponds to visual improvement
in the LSQ-DTV reconstructed images as shown in Fig.~\ref{fig:phgeommw}. Note that there is significant improvement
in the recovery of the muscle wall background between 500 and 2000 iterations of LSQ-DTV.

\subsection{Applying LSQ-DTV image reconstruction to high-resolution DBT imaging}
\label{sec:highres}

For low-resolution image reconstruction as presented in Sec.~\ref{sec:lowres}, executing an iterative image
reconstruction algorithm for hundreds or thousands of iterations may be practical because of the reduced
image resolution and grid size. Full resolution DBT, however, running LSQ-DTV to near convergence may not be feasible especially
for voxels with a low aspect ratio -- on the order of 100 microns in size.
The strategy proposed here
is to accurately solve the low-resolution LSQ-DTV problem image reconstruction, then to use the resulting
volume to initialize a gradient descent algorithm solving LSQ-Tik at high-resolution.

\begin{figure}
 \centering
        \includegraphics[width=0.45\textwidth]{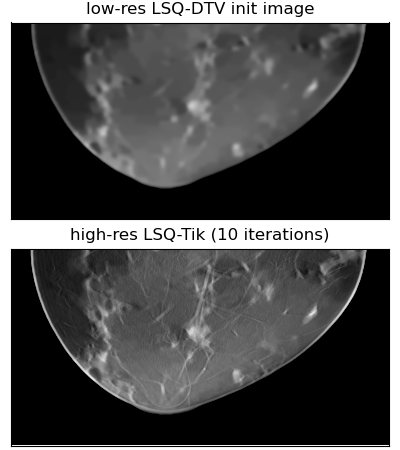}
        \includegraphics[width=0.45\textwidth]{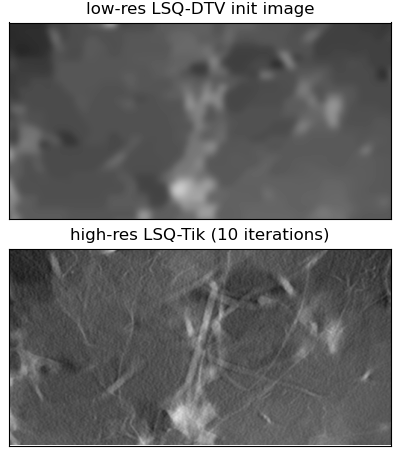}
 \caption{In-plane slice images of the initial low-resoluction LSQ-Tik image (Top row) and the high-resolution
 LSQ-Tik image after 10 iterations of gradient descent (Bottom row). The full slice is shown in the
 left column and a blow-up on a region-of-interest in the right column. The gray scale is
 $[0.45,0.9]$ cm$^{-1}$.}
\label{fig:gdvic}
\end{figure}

To demonstrate this approach, we use an anthropomorphic phantom generated by software from the Virtual
Imaging Clinical Trials for Regulatory Evaluate (VICTRE) project.\cite{VICTRE}
The low-resolution LSQ-DTV
DBT volume is reconstructed using the same parameter settings as the geometric-shapes phantom
study of Sec.~\ref{sec:lowres} except that the data-discrepancy constraint is tuned to $\epsilon= 0.003$.
The anthropomorphic phantom is smaller than the geometric-shapes phantom and the grid of 400 micron
voxels is 253$\times$132$\times$135, and the modeled detector has 1619$\times$591 pixels of 100 microns width.
As described in Eq.~(\ref{DTVopt3D}), the object model is a gaussian-blurred latent image $G[d] f$.
When creating the initial image for the high-resolution LSQ-Tik optimization, the latent image $f$
is upsampled to the high-resolution grid using nearest neighbor interpolation, then blurred with the gaussian kernel.
Note that, in order to avoid blockiness of the initial image, the upsampling needs to be performed first.
For the present example, we upsample to an in-plane resolution of 100 microns, and a depth resolution of 200 microns.
The initial image $h_0$ is specified on a grid of 1012$\times$528$\times$270 voxels. The gaussian blurring
used in forming $h_0$ is still specified by a width $d=(400,400,400)$ microns, which is the same value used
in the low-resolution LSQ-DTV computation.

To obtain the high-resolution LSQ-Tik reconstruction, gradient descent is performed
solving
\begin{equation}
\label{LSQTikHR}
\argmin_h \left\{ \tfrac{1}{2}\|R[c](XG[d]h-g)\|^2_2 + \alpha_\text{Tik} \tfrac{1}{2} \| G[d]h - h_0 \|^2_2  \right\},
\end{equation}
where $h_0$ is a prior image used for initialization and regularization; $G[d]$ symbolizes the gaussian used
to blur the latent image; and $\alpha_\text{Tik}$ is the regularization
parameter. In princple, if Eq.~(\ref{LSQTikHR}) is solved accurately, an initial image is not needed; but the intent
here is to run only a few gradient descent iterations in which case initializing with $h_0$ helps to obtain a useful image.
For comparison, high-resolution LSQ-Tik is solved with the LSQ-DTV initial image estimate (DTVinit-LSQ-Tik) and
with $h_0=0$ (zeroinit-LSQ-Tik).
For the results presented here, 10 iterations of gradient descent are performed for both DTVinit-LSQ-Tik
and zeroinit-LSQ-Tik.
Fig.~\ref{fig:gdvic} shows an in-plane
slice image for both the LSQ-DTV prior, $h_0$, and
DTVinit-LSQ-Tik image after 10 iterations of gradient descent, $h^{(k=10)}$. 
The gaussian blur width is set to $d=(100,100,200)$ microns and
the regularization parameter is set to $\alpha_\text{Tik} = 0.05$ based on noise-level in
the final image. From these images, it is clear that the low-resolution prior image is missing detailed
structures, which are recovered with DTVinit-LSQ-Tik on a high-resolution grid.

\begin{figure}
\centering
        \includegraphics[width=0.425\textwidth]{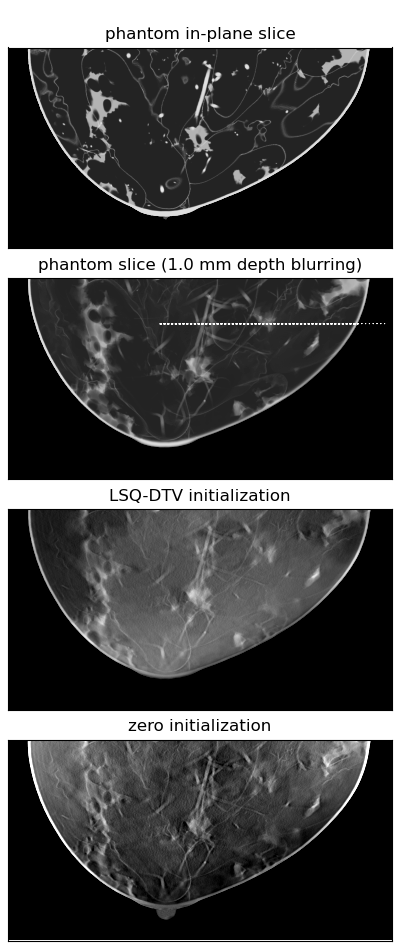}
        \includegraphics[width=0.425\textwidth]{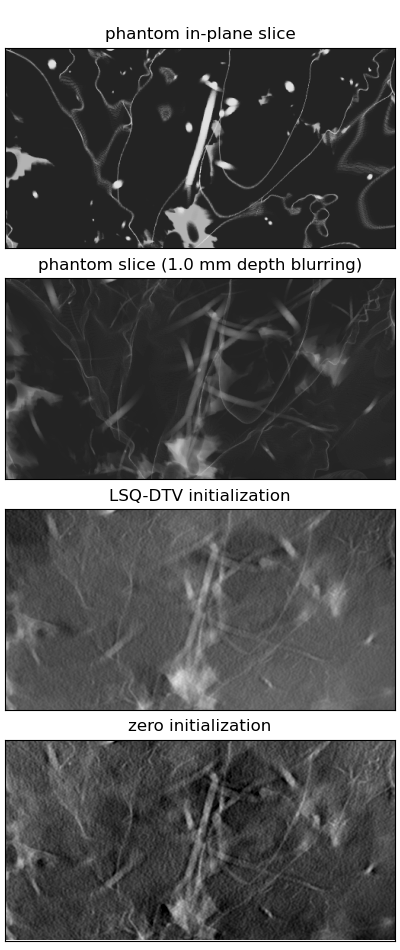}
\caption{Comparison of DTVinit-LSQ-Tik and zeroinit-LSQ-Tik.
The top row is an in-plane 100 micron slice of the anthropomorphic breast phantom. The second
row is the same phantom plane after depth blurring with a gaussian distribution of width 1 mm.
The third row, which is the same as the second row in Fig.~\ref{fig:gdvic}, shows the DTVinit-LSQ-Tik result and the
fourth row shows the same for zeroinit-LSQ-Tik.
The dotted line in the full image of the depth-blurred phantom image indicates the location of the profile
plotted in Fig.~\ref{fig:phvicprofile}.
The gray scale window is $[0.45,0.9]$ cm$^{-1}$ but the fourth row results have been scaled so that phantom
features are clearly visible in this gray scale window.}
\label{fig:phvic}
\end{figure}

\begin{figure}
\centering
        \includegraphics[width=0.45\textwidth]{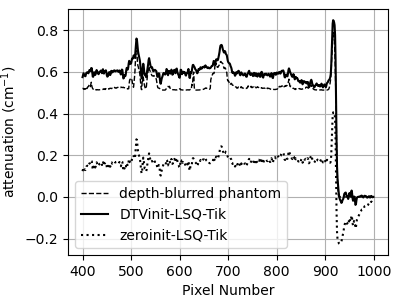}
\caption{Profile plot of high-resolution LSQ-Tik reconstructions shown in Fig.~\ref{fig:phvic}.
The location of the profile is indicated by the dashed line in Fig.~\ref{fig:phvic}.}
\label{fig:phvicprofile}
\end{figure}

The comparison of DTVinit-LSQ-Tik and zeroinit-LSQ-Tik  is shown in Fig.~\ref{fig:phvic}.
Along with this comparison are images of the phantom in a thin 100 micron slice and depth-blurred with a gaussian
of width 1 mm. For a more quantitative comparison, line profile curves of the comparison are shown in Fig.~\ref{fig:phvicprofile}.
The biggest differences that the LSQ-DTV prior makes is in the quantitative accuracy; DTVinit-LSQ-Tik
is within 10\% of the ground truth, while zeroinit-LSQ-Tik has the undershoot artifacts typical of
standard DBT image reconstruction and underestimates the ground truth by nearly a factor of 3.
The features
of the 1 mm depth-blurred phantom slice are seen in both DTVinit-LSQ-Tik and zeroinit-LSQ-Tik
reconstructed images. DTVinit-LSQ-Tik does
a better job than zeroinit-LSQ-Tik at suppressing structures that are in planes above and below the current viewing plane;
for example, the nipple is visible in the zeroinit-LSQ-Tik while it is not in the ground truth or DTVinit-LSQ-Tik images.
Furthermore, zeroinit-LSQ-Tik is more blotchy relative to DTVinit-LSQ-Tik due to out-of-focus features of structures displaced
in $z$ relative to the viewing plane.

\section{Discussion the DBT application}

The refined algorithmic framework allows for the development of multi-term convex optimization for imaging applications
by simplifying the selection of the myriad of parameters that control the step-size of the PDHG iteration.
This framework allows refinement of DTV-regularization for improving DBT imaging in terms
of quantitativeness and depth resolution as demonstrated
in Sec.~\ref{sec:3DDBT}. With this framework in hand, even more complex modeling becomes feasible
such as the use of dual-channel data fidelity optimization for further improving convergence in DBT
image reconstruction.\cite{Iyadomi2026} For the version of DTV-regularization, demonstrated in this work,
depth-blurring is controlled to the point where the resulting reconstruction is visually similar
to the ground truth with depth blurring of 1 mm for a wide-angle DBT setup. Quantitative accuracy is improved
but there is still a visible shift in the gray value with respect to the ground truth. Also, more empirical
studies are needed using anthropomorphic phantoms of different size and complexity.

\section{Conclusion}

Efficient PDHG implementation with matrix stacking for multi-term convex optimization is presented. The formalism
is demonstrated on an imaging application, namely image reconstruction in wide-angle DBT. This is the application
that motivated us to simplify the step parameter selection process. Specifically, the optimization
problem of interest, shown in Eq.~(\ref{DTVopt3D}), has 6 terms that each come with their own linear transform.
In this case, there are 13 parameters that specify the diagonal step-parameter matrices, $\Sigma$ and $T$, which
makes a grid search over these parameters infeasible. With the simplified parameter selection, the selection
of the 13 parameters is reduced to a determination of only 2 parameters. This simplification thus enables
the use of acceleration techniques that come with their own set of parameters; in our case we use the
predictor-corrector PDHG algorithm of He and Yuan,\cite{he2012convergence} which adds a third parameter.
 The step parameterization
scheme presented here is specific to optimization problems for X-ray tomography image reconstruction, where
the X-ray transform is dominant in terms of specifying the data-fidelity and in terms of computational effort.
Other tomographic applications may require a different scheme.

\ack{
This work is supported in part by NIH grant No. R01-CA287302.
The contents of this article are solely the responsibility of the authors and do not necessarily represent
the official views of the National Institutes of Health.
}

\providecommand{\newblock}{}

\end{document}